\documentclass[12pt]{article}
\usepackage[T2A]{fontenc}
\usepackage[cp1251]{inputenc}
\usepackage[russian]{babel}
\usepackage{amsmath,amssymb,amsthm,amsfonts,amscd}
\begin{document}

\centerline{\uppercase{\bf Parabolic polygons }}
\bigskip
\centerline{\bf Nilov Feodor }
\smallskip
\centerline{lyceum "The Second school"}

\bigskip 
{\bf Abstract.}
\footnote{This paper is prepared under the supervision of A.Zaslavski and is submitted 
to a prize of the Moscow Mathematical Conference of High-school Students. 
Readers are invited to send their remarks and reports to mmks@mccme.ru.}



                                          
\smallskip
{\bf Main Theorem.} {\it Two parabols have four common points. 
There exists a circle tangent to the sides of the obtained parabolic quadrilateral 
if and only if the diagonals of this quadrilateral are orthogonal. }

\smallskip
The proof of the Main Theorem is elementary and purely synthetic. 
It is based on the following lemma. 

\smallskip
{\bf Lemma 1.} {\it Assume that a parabola is tangent to a circle at points 
$A$ and $B$. 
A point $P$ of the plane lyes on the parabola if and only if 
the distance from the point $P$ to the line $AB$ equals to the length of the 
tangent from $P$ to the circle.}

\smallskip
{\bf Corollary 1.} {\it Each parabolic quadrilateral as above is affine equivalent 
to a parabolic quadrilateral which is simultaneously inscribed and 
circumscribed. }
 
\smallskip
{\bf Corollary 2.} {\it A parabola and a circle have four common points $A$, $B$, 
$C$, $D$. 
Denote by $L$ the intersection point of diagonals $AC$ and $BD$ 
of quadrilateral $ABCD$. 
Then the bissector of the angle $CLD$ is parallel to the axis of the parabola. }

\smallskip
{\bf Proposition 1.} {\it Two parabols have four common points lying on a 
circle. 
Then the intersection point of the axes of these parabols is the barycentre of 
the four points.}


\newpage

\centerline{\uppercase{\bf Параболические многоугольники }}
\bigskip
\centerline{\bf Нилов Фёдор }
\smallskip

\centerline{лицей "Вторая школа"}
\bigskip

В этой работе доказывается несколько теорем о криволинейном параболическом 
четырёхугольнике.

\smallskip
{\bf Основная Теорема.}
{\it Две параболы пересекаются в четырёх точках.
В полученный "параболический четырёхугольник" можно вписать окружность
тогда и только тогда, когда его диагонали перпендикулярны.  }


\smallskip
По-видимому, этот красивый факт неизвестен, что подтверждается мнением
авторов работ [1, 2].





\bigskip

\centerline{\bf Формулировки остальных результатов.}
\smallskip
\smallskip



{\it Параболическим четырёхугольником} 
назовём фигуру, являющуюся пересечением двух частей плоскости, каждая из 
которых ограничена параболой, причем эти две параболы пересекаются 
в четырех точек 

\smallskip
{\bf Следствие 1.}
{\it Любой параболический четырёхугольник можно перевести
аффинным преобразованием во вписанно-описанный параболический четырёхугольник.}



\smallskip
Из следствия 1 будет выведен еще один результат, следствие 4.
Но сначала сформулируем следующий интересный факт, который 
приводится к частному случаю следствия 4.

\smallskip
{\bf Следствие 2.}  {\it Парабола и окружность пересекаются
в четырёх точках $A$, $B$, $C$ и $D$. Обозначим точку пересечения диагоналей
четырёхугольника $ABCD$ через $L$. Тогда биссектриса угла $CLD$
параллельна оси параболы.}

\smallskip
Cледующее утверждение будет выведено из следствия 2.

\smallskip
{\bf Следствие 3.} {\it Две параболы пересекаются в четырех точках, 
лежащих на одной окружности. 
Тогда точка пересечения осей парабол совпадает с центром тяжести 
этих четырёх точек.}


Прямая $EF$ называется {\it осевой прямой} выпуклого четырёхугольника $ABCD$, 
если
она проходит через точку пересечения диагоналей четырёхугольника и пересекает
прямые, содержащие стороны $AB$ и $CD$ в точках $E$ и $F$, для которых
$AE / EB = FD / CF$. Осевая прямая четырёхугольника зависит от нумерации
вершин. 

Рассмотрим аффинное преобразование, которое переводит
четырёхугольник $ABCD$ во вписанный четырёхугольник
$A'B'C'D'$. Тогда осевая прямая четырёхугольника
$ABCD$ при таком преобразовании перейдёт в биссектрису угла $C'L'D'$, 
где $L'$ --- точка пересечения диагоналей $A'B'C'D'$.

\smallskip
{\bf Следствие 4.} {\it На параболе лежат четыре точки $A$, $B$, $C$ и $D$.
Тогда ось параболы параллельна осевой прямой четырёхугольника $ABCD$.}



\smallskip
Следующее утверждение будет выведено из леммы 1, сформулированной ниже.


\smallskip
{\bf Утверждение 1.} {\it Главные диагонали описанного параболического
шестиугольника пересекаются в одной точке.}


\smallskip
Следующее утверждение принадлежит А.В. Акопяну и А.А. Заславскому.
Я привожу другое доказательство, основанное на следствии 2 и лемме 1, 
сформулированной далее.

\smallskip
{\bf Теорема.} {\it Внутри окружности 
взята точка $X$. Через неё проводятся $N$ хорд, делящие плоскость на 
$2N$ равных углов. Через концы каждой хорды проводится парабола, 
касающаяся окружности в этих концах.
Тогда вершины параболического $2N$-угольника, 
 получающегося при пересечении
этих парабол, лежат на одной окружности.}






\bigskip
\centerline{\bf Доказательство части "только тогда" основной теоремы}
\smallskip
\smallskip





Сформулируем лемму 1, на которой основано доказательство основной теоремы.

\smallskip
{\bf Лемма 1.} {\it Дана парабола, касающаяся окружности в точках $A$ и $B$.
Произвольная точка $P$ плоскости лежит на этой параболе тогда и только тогда,
когда расстояние от точки $Р$ до прямой $AB$ равно длине касательной,
проведённой из точки $P$ к окружности.}

\smallskip
Доказательство этой леммы приведём позже.

\smallskip
{\it Доказательство части "только тогда" основной теоремы.}
Обозначим точки касания окружности, вписанной в параболический
четырёхугольник, с одной параболой через $K$ и $L$, а с другой ---
через $M$ и $N$.
Рассмотрим одну из вершин $A$ параболического четырёхугольника.
Из части "только тогда" леммы 1 следует, что расстояния от точки $A$ до прямых
$KL$ и $MN$ равны длине касательной, проведённой из $A$ к окружности.
Значит, вершинa $A$ равноудалена от прямых $KL$ и $MN$.
Аналогичное верно и для других вершин.
Поскольку прямые, содержащие биссектрисы углов,
образованных прямыми $KL$ и $MN$, являются диагоналями четырёхугольника $ABCD$,
то диагоналями этого четырёхугольника перпендикулярны.

\bigskip
\centerline{\bf Доказательство части "тогда" основной теоремы.}
\smallskip
\smallskip

\smallskip
Основная идея доказательства части 'тогда' основной теоремы 
принадлежит А.А. Заславскому.

\smallskip
{\it Доказательство части 'тогда' основной теоремы.} Пусть $ABCD$ --- параболический четырехугольник с перпендикулярными
диагоналями.
Обозначим через $L$ точку пересечения его диагоналей.
Известно, что проекции точки $L$ на прямые $AB$, $BC$, $CD$ и $DA$ лежат
на одной окружности. Докажем, что эта окружность вписана в параболический
четырёхугольник.
Можно считать, что существует такая прямая $l_1$, проходящая через точку $L$,
что расстояние $a$ от этой
прямой до
точки $A$ равно длине касательной $t_a$ из точки $A$ к окружности.

(Покажем, почему такая прямая существует. Будем считать, что центр $I$ этой окружности
лежит в той же полуплоскости относительно прямой $BD$, что и точка $A$.
Тогда $\angle ILA$ острый. Существование такой прямой по соображениям непрерывности следует из
$AL^2>t_a^2$. Для доказательства этого неравенства обозначим через $r$
радиус рассматриваемой окружности. Точка $L$
лежит внутри этой окружности. Поскольку $\angle ILA$ острый,
$AL^2>AI^2-IL^2>AI^2-r^2=t_a^2$.)

\smallskip


Обозначим расстояния от этой прямой до
каждой из точек $B$, $C$ и $D$ через $b$, $c$, $d$, а длины
касательных из точек $B$, $C$ и $D$ к окружности через $t_b$,
$t_c$ и $t_d$.


Докажем, что $t_b=b$.

Обозначим проекцию точки $L$ на прямую $AB$ через $H$.
Обозначим вторую точку пересечения окружности и прямой $AB$ через $X$.
Тогда
$$\displaystyle\frac{t_a^2}{AL^2}+\frac{t_b^2}{BL^2}=
\frac{AH\cdot AX}{AH\cdot AB}+\frac{BX\cdot BH}{BH\cdot AB}=
\frac{AX+XB}{AB}=1=$$
$$ =\sin^2\angle(l_1,AC)+\cos^2\angle(l_1,AC)=
\frac{a^2}{AL^2}+\frac{b^2}{BL^2}.$$
Это равенство выполнено, поскольку

(a) $LH$ --- высота прямоугольного
треугольника $ALB$, опущенная на гипотенузу $AB$, а значит
$AL^2=AH\cdot AB$ и $BL^2=BH\cdot AB$.

(b) $AH\cdot AX=t_a^2$ и $BX\cdot BH=t_b^2$.

Поэтому $a=t_a$ влечет
$b=t_b$.
Аналогично, $c=t_c$ и $d=t_d$.


Поскольку для прямой $l_1$ выполнено $t_a=a$, $t_b=b$, $t_c=c$ и $t_d=d$, то
то прямая $l_2$, симметричная ей относительно $AC$
тоже обладает этим же свойством.
Поэтому из части "тогда" леммы 1 следует,
что точки $A$, $B$, $C$ и $D$ лежат на двух параболах,
касающихся окружности в точках пересечения прямых $l_1$ и $l_2$ с окружностью.
Поэтому в параболический четырёхугольник $ABCD$ можно вписать окружность.


\bigskip
\centerline{\bf Лемма 2.}
\smallskip
\smallskip

Докажем лемму, на которой основано доказательство
части "только тогда" леммы 1.

\smallskip

{\bf Лемма 2.} {\it По параболе движется точка $P$.
Пусть $AB$ --- хорда параболы,
параллельная её директрисе.

(a) точка $C$ пересечения перпендикуляров, восстановленных в
точках $A$ и $B$ к прямым $PA$ и $PB$, движется по прямой, параллельной $AB$.


Обозначим проекцию
точки $Р$ на прямую $АВ$ через $H$, а проекции точки $H$
на прямые $AP$ и $BP$ через $K$ и $L$. Точки $A$, $B$, $K$ и $L$
лежат на одной окружности $\omega$, поскольку
$\angle KLP = \angle KHP = \angle HAP$. Обозначим через $O$ центр окружности
$\omega$. Тогда

(b) окружность $\omega$
не зависит от положения точки $Р$.

(c) окружность $\omega$ касается параболы в точках $A$ и $B$.}

\smallskip

Для доказательства леммы 2 (a) приведём определения и
сформулируем известную лемму.

\smallskip

Пучком $[A]$ прямых называется множество прямых, проходящих через точку $A$.

\smallskip


Соответствие между двумя пучками прямых называется
{\it проективным}, если двойное отношение любых четырёх прямых из одного пучка
равно двойному отношению четырёх соответствующих прямых из другого пучка.

\smallskip

{\bf Лемма Соллертинского.} {\it Дано проективное соответствие между пучками
прямых $[А]$ и $[В]$. Тогда все точки $P$,
являющиеся пересечением соответствующих прямых
из пучков $[A]$ и $[B]$, принадлежат коническому сечению
(возможно, вырожденному).}

\smallskip

{\it Доказательство леммы 2 (a).}
Рассмотрим прямую, проходящую через точку $C$ и
параллельную прямой $AB$. Пусть точка $C'$ движется по этой прямой.
Обозначим точку пересечения перпендикуляров, восстановленных в точках
$A$ и $B$ к прямым $C'A$ и $C'B$, через $P'$.
Очевидно, что соответствие $AC'\to BC'$ между прямыми из пучков $[A]$ и $[B]$
является проективным.
Соответствие $AC'\to AP'$ между прямыми из пучков $[A]$
является проективным, поскольку угол между любыми двумя соответствующими
прямыми прямой.
Аналогично, соответствие $BC'\to BP'$ между прямыми из пучков $[B]$
является проективным.
Значит, соответствие $AP'\to BP'$ между прямыми из пучков $[A]$ и $[B]$
является проективным.

Поэтому из леммы Соллертинского следует,
что точки $P'$ лежат на конике (возможно, вырожденной). 
Легко убедится 
в том, что кривая $\gamma$,
по которой движется точка $P'$, не является вырожденной коникой. 
Для этого 
достаточно рассмотреть пять точек кривой $\gamma$, никакие три из которых 
не лежат на одной прямой. Это точки $A$, $B$, две точки пересечения серединного 
перпендикуляра к отрезку $AB$ c кривой $\gamma$ (одна из которых конечна, а 
другая бесконечно удалена) и произвольная пятая точка кривой $\gamma$. 


Заметим, что кривой $\gamma$ принадлежит ровно одна бесконечно удалённая точка
(эта точка является пересечением перпендикуляров, восстановленных в точках
$A$ и $B$ к отрезку $AB$). Коника, которой принадлежит ровно одна
бесконечно удалённая точка, является параболой.
Значит, $\gamma$ является параболой, ось которой перпендикулярна прямой
$AB$. Поскольку эта парабола проходит через точки $A$, $B$ и $P$ и её ось
параллельна оси данной параболы, эта парабола совпадает с данной.
Поэтому точка $C$ движется по прямой, параллельной $AB$.

\smallskip

{\it Доказательство леммы 2 (b).}
Рассмотрим случай, когда точка $P$ находится в
той же полуплоскости относительно прямой $AB$, что и
вершина параболы (cлучай, когда точка $P$ находится в другой
полуплоскости относительно прямой $AB$ доказывается аналогично).
Серединные перпендикуляры к отрезкам $AK$ и $LB$
являются средними линиями трапеций $AKHC$ и $BLHC$.
Следовательно, они пересекаются в середине отрезка
$HC$.


Значит, середина отрезка $HC$ совпадает с точкой $O$.
Обозначим проецию
точки $C$ на прямую $AB$ через $S$. Из леммы 2 (a) следует,
что длина отрезка $CS$ не зависит от положения точки $P$
на параболе. 

Обозначим середину отрезка $AB$ через $M$.
Точка $O$ находится на серединном
перпендикуляре к отрезку $AB$, причём $OM = 1/2 CS$.
Следовательно, положение точки $O$ не зависит от выбора
точки $P$ на параболе.

\smallskip
{\it Доказательство леммы 2 (c).}
Рассмотрим точку $P = A$. Тогда прямая $AP$ будет
касательной к параболе в точке $A$. Из того, что
прямые $AO$ и $AP$ перпендикулярны и точка $O$ лежит на
серединном перпендикуляре к отрезку $AB$ следует, что точка $O$ совпадёт
с центром окружности, касающейся параболы в точках
$A$ и $B$.

Из леммы 2 (b) и вышесказанного
следует, что для произвольного положения
точки $P$ на параболе, точка $O$ является
центром окружности, касающейся параболы в точках $A$ и $B$.
Значит, окружность $\omega$
совпадает с окружностью, касающейся параболы
в точках $A$ и $B$.

\bigskip

\centerline{\bf Доказательство леммы 1.}
\smallskip
\smallskip

{\it Доказательство части "только тогда" леммы 1.}
Обозначим проекцию
точки $Р$ на прямую $АВ$ через $H$, а проекции точки $H$
на прямые $AP$ и $BP$ через $K$ и $L$.
Поскольку треугольники $HPB$ и $HPL$ подобны и
окружность, касающаяся параболы в точках $A$ и $B$,
проходит через точки $K$ и $L$ (это следует из пункта (c) леммы 2),
то $PH=\sqrt {PL \cdot PB}$. Значит, отрезок $PH$ равен длине касательной,
проведённой из точки $P$ к окружности, касающейся параболы в точках $A$ и $B$.

\smallskip
{\it Доказательство части "тогда" леммы 1.}
Предположим противное.
Пусть точка $P$ не лежит на параболе, касающейся окружности в точках
$A$ и $B$, но расстояние от точки $P$ до
прямой $AB$ равно длине касательной, проведённой из точки $P$ к окружности.
Рассмотрим случай, когда точка $P$ находится в той же полуплоскости
относительно прямой $AB$, что и вершина параболы.
Обозначим проекцию точки $P$ на прямую $AB$ через $H$. Обозначим точку
пересечения прямой $PH$ с параболой через $P'$. Точки $P$ и $P'$ различны,
поскольку точка $P$ не лежит на параболе. Пусть $PX$ и $P'X'$ --- отрезки
касательных, проведённых из точек $P$ и $P'$ к окружности.
Тогда $PH=PX$ по условию. Из части "только тогда" леммы 1 следует, что
$P'H=P'X'$. Обозначим центр и радиус окружности через $I$ и $r$.
Обозначим проекцию точки $I$ на прямую $PH$ через $H'$.
Тогда
$$\displaystyle PH^2-P'H^2=PX^2-P'X'^2=(PI^2-r^2)-(P'I^2-r^2)=PI^2-P'I^2=$$
$$=PH'^2-P'H'^2 = (PH+HH')^2-(P'H+HH')^2=$$
$$=PH^2-P'H^2+2HH'\cdot(PH-P'H).$$

Следовательно, $2HH'\cdot(PH-P'H)=0$. Но $PH\ne P'H$ и  $HH'\ne 0$, поскольку хорда
$AB$ не является диаметром. Противоречие. Случай, когда точка $P$ находится
в другой полуплоскости относительно прямой $AB$, доказывается аналогично.

\bigskip
\centerline{\bf Доказательства остальных результатов.}
\smallskip
\smallskip


{\it Доказательство следствия 1.} 
Сформулируем следующий известный факт.

\smallskip
{\bf Теорема.} {\it Параболический четырёхугольник является вписанным тогда и
только тогда, когда оси образующих его парабол перпендикулярны.}

\smallskip
Аффинным преобразованием переведём
оси парабол и диагонали параболического четырёхугольника
в две пары перпендикулярных прямых. Как известно,
при произвольном аффинном преобразовании образ оси параболы
параллелен оси образа параболы. Значит,
оси образа параболического четырёхугольника перпендикулярны.
Поэтому из вышеприведённой теоремы и основной теоремы следует, что параболический
четырёхугольник при таком преобразовании переходит во
вписанно-описанный.

\smallskip
{\it Доказательство следствия 2.} Обозначим точки пересечения биссектрисы
угла $ALD$ с прямыми $BC$ и $AD$ через $M$ и $N$. Тогда
$BM / MC = LB / LC = AL / LD = AN / ND$. Значит,
биссектрисы углов, образованных диагоналями
вписанного четырёхугольника являются осевыми
прямыми этого четырёхугольника. Поэтому из следствия 4 (доказательство
см. ниже) вытекает следствие 2.

\smallskip
{\it Доказательство следствия 3.} Рассмотрим систему координат в которой одна
из парабол
имеет уравнение
$y = kx^2$. Обозначим координаты
точек $A$, $B$, $C$ и $D$ в этой системе через $(x_1, y_1)$, $(x_2, y_2)$,
$(x_3, y_3)$ и $(x_4, y_4)$. Тогда $y_1 = kx_1^2$, $y_2 = kx_2^2$,
$y_3 = kx_3^2$, $y_4 = kx_4^2$. Обозначим коэффициенты $k$ прямых
$AB$ и $CD$ через $k_1$
и $k_2$. По следствию 2 $k_1 = -k_2$. Поэтому                                                                                                           
$$\displaystyle \frac{y_2 - y_1} {x_2 - x_1} = \frac {k (x_1 + x_2)} {2} = \frac {- (y_3 - y_4)} {x_3 - x_4} = \frac {- k (x_3 + x_4)} {2}.$$
Заметим,
что $\displaystyle \frac {k (x_1 + x_2)} {2}$ и $\displaystyle \frac {- k (x_3 + x_4)} {2}$ являются
абциссами середин отрезков $AB$ и $CD$ в заданной системе координат.
Отсюда следует, что ось выбранной параболы проходит через середину отрезка,
соединяющего середины отрезков $AB$ и $CD$. Поэтому она проходит через центр
тяжести четырёхугольника $ABCD$. То же самое можно сказать и про ось
другой параболы. Значит, оси парабол пересекаются в центре тяжести
четырёхугольника $ABCD$.

\smallskip
{\it Доказательство следствия 4.} Рассмотрим параболический четырёхугольник с вершинами
в точках $A$, $B$, $C$ и $D$. Используя следствие 1, переведём соответствующий
параболический четырёхугольник во вписанно-описанный.
Его осевая прямая параллельна оси параболы, являющейся образом исходной параболы.
Следовательно, и осевая прямая
четырёхугольника $ABCD$ параллельна оси параболы.












\smallskip
{\it Доказательство утверждения 1.} Обозначим точки касания окружности с одной
параболой через $A_1$ и $A_2$, со второй --- через $B_1$ и $B_2$,
с третьей --- через $C_1$ и $C_2$. Из леммы 1
следует, что главные диагонали полного описанного
параболического шестиугольника являются биссектрисами
треугольника, образованного прямыми $A_1A_2$, $B_1B_2$ и $C_1C_2$. Значит,
они пересекаются в одной точке.

\smallskip
{\it Доказательство теоремы Акопяна и Заславского.} Докажем эту теорему для $N = 3$.
Из доказательства будет видно, что общий
случай доказывается аналогично.

Следующая лемма очевидна.

\smallskip
{\bf Лемма.} {\it Через точку внутри окружности проведены две хорды.
Тогда парабола, касающиеся окружности в концах первой хорды, и парабола,
касающаяся окружности в концах второй хорды, пересекаются в четырёх точках.}

\smallskip
Из этой леммы следует,
что главные диагонали шестиугольника являются
биссектрисами углов, на которые хорды делят плоскость.
Поэтому хорды и главные диагонали делят плоскость
на 12 равных углов. Рассмотрим две главные диагонали
шестиугольника, концы которых принадлежат одной параболе.
Тогда хорда, принадлежащая этой параболе,
является биссектрисой угла, образованного этими диагоналями.
Из этого и следствия 2 следует, что концы
этих главных диагоналей, являющиеся вершинами шестиугольника,
лежат на окружности. Мы доказали, что концы любых двух главных диагоналей
шестиугольника лежат на одной окружности. Поэтому все вершины
шестиугольника лежат на одной окружности.



\smallskip

{\it Покажем, как можно построить с помощью циркуля и линейки
бесконечно много точек
параболического четырёхугольника с
вершинами в четырёх данных точках.} Известно, что через четыре
точки можно провести только две параболы.
Приведём план этого построения.
Обозначим через $A$, $B$, $C$ и $D$ четыре данные точки.
Построим две осевые прямые четырёхугольника $ABCD$.
Они будут параллельны осям парабол, образующих параболический
четырёхугольник. Обозначим эти параболы через $p_1$ и $p_2$.
Проведём окружность через точки $A$, $B$ и $C$.
Из следствия 2 следует, что если эта окружность и
парабола $p1$ (или $p2$) пересекаются в четырёх точках, то
прямые, соединяющие эти точки, образуют равные углы с
осевыми прямыми четырёхугольника $ABCD$, следовательно,
можно построить четвёртую точку пересечения
окружности с параболой $p1$ (или $p2$). Таким образом, построены пять
точек параболы $p1$ (или $p2$). Значит, можно построить бесконечно
много точек, принадлежащих параболическому четырёхугольнику.

\smallskip
{\bf Благодарность.}
\smallskip

Автор благодарит А.А. Заславского за конструктивное обсуждение работы
и А.Б. Скопенкова за ценные замечания при подготовке текста.

\smallskip


[1] А.В. Акопян, А.А. Заславский.
{\it Геометрические свойства кривых второго порядка.} М.: МЦНМО, 2007.
English translation is in preparation.

[2] С.В. Маркелов. 
{\it Парабола как окружность.} // Десятая летняя конференция 
Турнира городов. М.: МЦНМО, 1999. С. 36-42, 112-123.

\end{document}